\documentclass[reqno,onecolumn,12pt]{amsart}
\usepackage{amsfonts}
\usepackage{amssymb}
\usepackage{amsmath}
\usepackage{graphicx}
\usepackage{amscd,amsthm}

\newtheorem{theorem}{Theorem}
\theoremstyle{plain}

\newtheorem{corollary}{Corollary}

\newtheorem{definition}{Definition}

\newtheorem{remark}{Remark}

\numberwithin{equation}{section}

\input{tcilatex}

\begin{document}
\author{}
\title{}
\maketitle

\begin{center}
\pagestyle{myheadings} \thispagestyle{empty} 
\markboth{\bf Serkan Araci, Mehmet Acikgoz, Hassan Jolany and Armen Bagdasaryan}
{\bf On the properties of q-Bernstein-type polynomials}

\textbf{\large On The Properties Of }\emph{q}\textbf{\large -Bernstein-Type
Polynomials}\bigskip

\textbf{Serkan Araci$^{1,\ast }$, Mehmet Acikgoz$^{2}$, Hassan Jolany$^{3}$}%
\bigskip \textbf{\ and Armen Bagdasaryan}$^{4,5}$

$^{1}$Faculty of Economics, Administrative and Social Sciences, Hasan
Kalyoncu University, 27410 Gaziantep, Turkey

\textbf{E-Mail: mtsrkn@hotmail.com}\\[0pt]

\textbf{$^{\ast }$Corresponding Author}\\[2mm]

$^{2}$Department of Mathematics, Faculty of Science and Arts, University of
Gaziantep, 27310 Gaziantep, Turkey

\textbf{E-Mail: acikgoz@gantep.edu.tr}

$^{3}$Universit\'{e} des Sciences et Technologies de Lille UFR de Math\'{e}%
matiques Laboratoire Paul Painlev\'{e} CNRS-UMR 8524 59655 Villeneuve d'Ascq
Cedex/France

\textbf{E-Mail: hassan.jolany@math.univ-lille1.fr}

$^{4}$Russian Academy of Sciences, Institute for Control Sciences
Profsoyuznaya 65, 117997 Moscow, Russia

$^{5}$American University of the Middle East, College of Engineering,
Department of Applied Mathematics Kuwait City, Block 3, Egaila, Kuwait

\textbf{E-Mail: bagdasar@member.ams.org}

\hspace{0.5cm}

\textbf{\large Abstract}
\end{center}

\begin{quotation}
The aim of this paper is to give a new approach to modified $q$-Bernstein
polynomials for functions of several variables. By using these polynomials,
the recurrence formulas and some new interesting identities related to the
second Stirling numbers and generalized Bernoulli polynomials are derived.
Moreover, the generating function, interpolation function of these
polynomials of several variables and also the derivatives of these
polynomials and their generating function are given. Finally, we get new
interesting identities of modified $q$-Bernoulli numbers\ and $q$-Euler
numbers applying $p$-adic $q$-integral representation on $%
\mathbb{Z}
_{p}$ and $p$-adic fermionic $q$-invariant integral on $%
\mathbb{Z}
_{p}$, respectively, to the inverse of $q$-Bernstein polynomials.
\end{quotation}



\section{\textbf{Introduction, Definitions and Notations}}

The Bernstein polynomials, named after their creater S. N. Bernstein in
1912, have been studied by many researchers for a long time. Recently
Acikgoz and Araci have originally defined the generating function of
Bernstein polynomials and analysed their interesting properties arising from
that generating function, and also the generating function of Bernstein
polynomials in two dimensional are defined by the same authors (see \cite%
{Acikgoz1}, \cite{Acikgoz2}, \cite{Acikgoz3}). Next, Simsek and Acikgoz have
constructed a generating function of ($q$-) Bernstein type polynomials based
on the $q$-analysis, \cite{SimsekAcikgoz1}, and gave some new relations
related to these polynomials, Hermite polynomials, Bernoulli polynomials of
higher order and the second kind Stirling numbers. Interpolation function of
($q$-) Bernstein type polynomials is defined by applying Mellin
transformation to this generating function. In \cite{Kim 1}, Kim-Choi-Kim
have studied on the $k$-dimensional generalization of $q$-Bernstein
polynomials, in which they have given some interesting properties of the $k$%
-dimensional generalization of $q$-Bernstein polynomials (see\cite{Kim 1}).
Our generalization of $q$-Bernstein polynomials are different from the $k$%
-dimensional generalization of $q$-Bernstein polynomials of Kim-Choi-Kim. In
the present paper, we also derived some interesting properties of our
generalization of $q$-Bernstein polynomials. Recent works including integral
representations and properties of Stirling numbers of the first kind \cite%
{Qi}, formulae for the $q$-Bernstein polynomials and $q$-deformed binomial
distributions \cite{Kim5}, integral representations for the Gamma function,
the Beta Function, and the double Gamma function \cite{Choi 4}, irregular
prime power divisors of the Bernoulli numbers \cite{Kellner}, application of
a composition of generating functions for obtaining explicit formulas of
polynomials \cite{Kruchinin}, hyperharmonic series involving Hurwitz zeta
function \cite{Mezo}, $p$-adic $q$-deformed fermionic integrals in the $p$%
-adic integer ring \cite{Araci} have been investigated extensively.

We are now in a position to give some definitions and some properties of
Bernstein polynomials of several variables with their generating function.

Let $C\left( \mathcal{D}^{w}\right) $ denotes the set of continuous
functions on $\mathcal{D}^{w}$, in which $\mathcal{D}^{w}$ and $\mathcal{D}$
mean $\underset{w-\text{times}}{\underbrace{\mathcal{D}\times \mathcal{D}%
\times ...\times \mathcal{D}}}$ and $\left[ 0,1\right] $, respectively. For $%
f\in C\left( \mathcal{D}^{w}\right) $, we have%
\begin{eqnarray}
\mathcal{B}_{n_{1},n_{2},\cdots ,n_{w}}\left( f;x_{1},x_{2},\cdots
,x_{w}\right)
&:&=\dsum\limits_{k_{1}=0}^{n_{1}}\dsum\limits_{k_{2}=0}^{n_{2}}\cdots
\dsum\limits_{k_{w}=0}^{n_{w}}f\left( \frac{k_{1}}{n_{1}},\frac{k_{2}}{n_{2}}%
,\cdots ,\frac{k_{w}}{n_{w}}\right)  \notag \\
&&\times B_{k_{1},k_{2},\cdots ,k_{w};n_{1},n_{2},\cdots ,n_{w}}\left(
x_{1},x_{2},\cdots ,x_{w}\right)  \label{Equation1}
\end{eqnarray}%
where $\mathcal{B}_{n_{1},n_{2},\cdots ,n_{w}}\left( f;x_{1},x_{2},\cdots
,x_{w}\right) $ is called the Bernstein operator of several variables of
order $\dsum\limits_{i=1}^{w}n_{i}$ for $f$. For $k_{i},n_{i}\in 
\mathbb{N}
_{0}$ with $i=1,2,\cdots ,w,$ the Bernstein polynomials of several variables
of degree $\dsum\limits_{i=1}^{w}n_{i}$ is defined by%
\begin{equation}
B_{k_{1},k_{2},\cdots ,k_{w};n_{1},n_{2},\cdots ,n_{w}}\left(
x_{1},x_{2},\cdots ,x_{w}\right) =\dprod\limits_{i=1}^{w}\left( \binom{n_{i}%
}{k_{i}}x_{i}^{k_{i}}\left( 1-x_{i}\right) ^{n_{i}-k_{i}}\right) \text{,}
\label{Equation2}
\end{equation}%
where $\binom{n}{k}=\frac{n\left( n-1\right) \cdots \left( n-k+1\right) }{k!}
$ and $x_{i}\in \mathcal{D}$ for $i=1,2,...,w$. These polynomials satisfy
the following relation%
\begin{equation*}
B_{k_{1},k_{2},\cdots ,k_{w};n_{1},n_{2},\cdots ,n_{w}}\left(
x_{1},x_{2},\cdots ,x_{w}\right)
=\dprod\limits_{i=1}^{w}B_{k_{i},n_{i}}\left( x_{i}\right)
\end{equation*}%
and they have form a partition of unity; that is:%
\begin{equation*}
\dsum\limits_{k_{1}=0}^{n_{1}}\dsum\limits_{k_{2}=0}^{n_{2}}\cdots
\dsum\limits_{k_{w}=0}^{n_{w}}B_{k_{1},k_{2},\cdots
,k_{w};n_{1},n_{2},\cdots ,n_{w}}\left( x_{1},x_{2},\cdots ,x_{w}\right) =1.
\end{equation*}

By using the definition of Bernstein polynomials for functions of several
variables, it is not difficult to prove the property given above as 
\begin{equation*}
\dsum\limits_{k_{1}=0}^{n_{1}}\dsum\limits_{k_{2}=0}^{n_{2}}\cdots
\dsum\limits_{k_{w}=0}^{n_{w}}\dprod\limits_{i=1}^{w}B_{k_{i},n_{i}}\left(
x_{i}\right) =1.
\end{equation*}

Also, $B_{k_{1},k_{2},\cdots ,k_{w};n_{1},n_{2},\cdots ,n_{w}}\left(
x_{1},x_{2},\cdots ,x_{w}\right) =0$ for $k_{i}>n_{i}$ with $i=1,2,...,w$,
because $\binom{n_{i}}{k_{i}}=0.$ There are $\dprod\limits_{i=1}^{w}\left(
n_{i}+1\right) $, $\dsum\limits_{i=1}^{w}n_{i}$-th degree Bernstein
polynomials.

Many researchers have studied the Bernstein polynomials of two variables in
approximation theory (see \cite{Buyukyazici1}, \cite{Buyukyazici2}). But
nothing was known about the generating function of these polynomials. Note
that for $k_{i},n_{i}\in 
\mathbb{N}
_{0}$\ and\textbf{\ }$x_{i}\in \mathcal{D}$ with $i=1,2,...,w$,\textbf{\ }we
obtain the generating function for $B_{k_{1},k_{2},\cdots
,k_{w};n_{1},n_{2},\cdots ,n_{w}}\left( x_{1},x_{2},\cdots ,x_{w}\right) $
as follows:%
\begin{gather}
F_{k_{1},k_{2},\cdots ,k_{w}}\left( t;x_{1},x_{2},\cdots ,x_{w}\right)
=\dprod\limits_{i=1}^{w}\frac{\left( tx_{i}\right) ^{k_{i}}}{k_{i}!}%
e^{wt-t\dsum\limits_{i=1}^{w}x_{i}}  \label{Equation2500} \\
=\dsum\limits_{n_{1}=k_{1}}^{\infty }\dsum\limits_{n_{2}=k_{2}}^{\infty
}\times \cdots \times  \notag \\
\dsum\limits_{n_{w}=k_{w}}^{\infty }B_{k_{1},k_{2},\cdots
,k_{w};n_{1},n_{2},\cdots ,n_{w}}\left( x_{1},x_{2},\cdots ,x_{w}\right)
\dprod\limits_{i=1}^{w}\frac{t^{n_{i}}}{n_{i}!}  \notag
\end{gather}%
where

\begin{equation*}
B_{k_{1},k_{2},\cdots ,k_{w};n_{1},n_{2},\cdots ,n_{w}}\left(
x_{1},x_{2},\cdots ,x_{w}\right) =\left\{ 
\begin{array}{ccc}
\dprod\limits_{i=1}^{w}\binom{n_{i}}{k_{i}}x_{i}^{k_{i}}\left(
1-x_{i}\right) ^{n_{i}-k_{i}} & \text{if} & n_{i}\geq k_{i}\text{,} \\ 
0 & \text{if} & n_{i}<k_{i}\text{,}%
\end{array}%
\right.
\end{equation*}%
for $k_{i}\in 
\mathbb{N}
_{0}$ and $x_{i}\in \mathcal{D}$, for $i=1,2,...,w$.

\begin{remark}
By substituting $w=1$ into (\ref{Equation2500}), we get a special case of $%
F_{k_{1},k_{2},\cdots ,k_{w}}\left( t;x_{1},x_{2},\cdots ,x_{w}\right) $
which was proved by Acikgoz and Araci (for details, see\ \cite{Acikgoz1})%
\begin{equation*}
F_{k_{1}}\left( t,x_{1}\right) =\frac{\left( tx_{1}\right) ^{k_{1}}e^{t}}{%
k_{1}!e^{tx_{1}}}=\sum_{n_{1}=k_{1}}^{\infty }B_{k_{1},n_{1}}\left(
x_{1}\right) \frac{t^{n_{1}}}{n_{1}!}.
\end{equation*}
\end{remark}

Let $0<q<1$. Define the $q$-number of $x$ by $[x]_{q}:=\frac{1-q^{x}}{1-q}$
\ and \ $[x]_{-q}:=\frac{1-\left( -q\right) ^{x}}{1+q},$ (see \cite{Araci1},%
\cite{Araci4},\cite{Araci5},\cite{Kim-Lee-Chae Jang},\cite{Kim 1},\cite{Kim
2},\cite{Kim 3},\cite{Kim 4},\cite{Ostrovska},\cite{Phillips},\cite%
{SimsekAcikgoz1} for details and related facts). Note that $\underset{%
q\rightarrow 1^{-}}{\lim }[x]_{q}=x$. \cite{Kim-Lee-Chae Jang} is actually
motivated the authors to write this paper and they have extended all results
given in \cite{Kim-Lee-Chae Jang} to modified $q$-Bernstein polynomials of
several variables.

\section{\textbf{The Modified }$q$\textbf{-Bernstein Polynomials for
Functions of Several Variables}}

For $0<q<1$, we consider 
\begin{gather*}
F_{k_{1},k_{2},\cdots ,k_{w}}\left( t,q;x_{1},x_{2},\cdots ,x_{w}\right)
=\dprod\limits_{i=1}^{w}\frac{\left( t\left[ x\right] _{q}\right) ^{k_{i}}}{%
k_{i}!}e^{t\dsum\limits_{i=1}^{w}\left[ 1-x_{i}\right] _{q}} \\
=\dsum\limits_{n_{1}=k_{1}}^{\infty }\dsum\limits_{n_{2}=k_{2}}^{\infty
}\times \cdots \times \\
\dsum\limits_{n_{w}=k_{w}}^{\infty }B_{k_{1},k_{2},\cdots
,k_{w};n_{1},n_{2},\cdots ,n_{w}}\left( x_{1},x_{2},\cdots ,x_{w};q\right)
\dprod\limits_{i=1}^{w}\frac{t^{n_{i}}}{n_{i}!}
\end{gather*}%
where $k_{i},n_{i}\in 
\mathbb{N}
_{0}$ and $x_{i}\in \mathcal{D}$ for $i=1,2,...,w$. We note that 
\begin{equation*}
\underset{q\rightarrow 1^{-}}{\lim }F_{k_{1},k_{2},\cdots ,k_{w}}\left(
t,q;x_{1},x_{2},...,x_{w}\right) =F_{k_{1},k_{2},\cdots ,k_{w}}\left(
t;x_{1},x_{2},...,x_{w}\right) .
\end{equation*}

\begin{definition}
We define the generating function of modified $q$-Bernstein polynomials for
functions of several variables as follows: 
\begin{gather}
F_{k_{1},k_{2},\cdots ,k_{w}}\left( t,q;x_{1},x_{2},\cdots ,x_{w}\right)
=\dprod\limits_{i=1}^{w}\frac{\left( t\left[ x_{i}\right] _{q}\right)
^{k_{i}}}{k_{i}!}e^{t\dsum\limits_{i=1}^{w}\left[ 1-x_{i}\right] _{q}}
\label{Equation58} \\
=\dsum\limits_{n_{1}=k_{1}}^{\infty }\dsum\limits_{n_{2}=k_{2}}^{\infty
}\times \cdots \times  \notag \\
\dsum\limits_{n_{w}=k_{w}}^{\infty }B_{k_{1},k_{2},\cdots
,k_{w};n_{1},n_{2},\cdots ,n_{w}}\left( x_{1},x_{2},\cdots ,x_{w};q\right)
\dprod\limits_{i=1}^{w}\frac{t^{n_{i}}}{n_{i}!}  \notag
\end{gather}%
where $k_{i},n_{i}\in 
\mathbb{N}
_{0}$ and $x_{i}\in \mathcal{D}$ with $i=1,2,...,w.$
\end{definition}

By using Taylor expansion of $e^{t\dsum\limits_{i=1}^{w}\left[ 1-x_{i}\right]
_{q}}$and the comparing coefficients on the both sides~in (\ref{Equation58}%
), we get the following Corollary.

\begin{corollary}
\label{Corollary 2.2.}For $k_{i},n_{i}\in 
\mathbb{N}
_{0}$ and $x_{i}\in \mathcal{D}$ for $i=1,2,...,w$, we have%
\begin{gather}
B_{k_{1},k_{2},\cdots ,k_{w};n_{1},n_{2},\cdots ,n_{w}}\left(
x_{1},x_{2},\cdots ,x_{w};q\right)  \notag \\
=\left\{ 
\begin{array}{ccc}
\dprod\limits_{i=1}^{w}\binom{n_{i}}{k_{i}}%
[x_{i}]_{q}^{k_{i}}[1-x_{i}]_{q}^{n_{i}-k_{i}} & \text{if} & n_{i}\geq k_{i}%
\text{,} \\ 
0 & \text{if} & n_{i}<k_{i}\text{.}%
\end{array}%
\right. .  \label{Equation59}
\end{gather}
\end{corollary}

\begin{theorem}
\textbf{(Recurrence Formula for }$B_{k_{1},k_{2},\cdots
,k_{w};n_{1},n_{2},\cdots ,n_{w}}\left( x_{1},x_{2},\cdots ,x_{w};q\right) $%
\textbf{)} For $k_{i},n_{i}\in 
\mathbb{N}
_{0},\ x_{i}\in \mathcal{D}$ and $i=1,2,...,w$, we have%
\begin{equation}
B_{k_{1},k_{2},\cdots ,k_{w};n_{1},n_{2},\cdots ,n_{w}}\left(
x_{1},x_{2},\cdots ,x_{w};q\right) =\dprod\limits_{i=1}^{w}\left(
[1-x_{i}]_{q}B_{k_{i};n_{i}-1}\left( x_{i};q\right)
+[x_{i}]_{q}B_{k_{i}-1;n_{i}-1}\left( x_{i};q\right) \right) \text{.}
\label{Equation 60}
\end{equation}

\begin{proof}
By using the definition of Bernstein polynomials for functions of several
variables, we have 
\begin{eqnarray*}
B_{k_{1},k_{2},\cdots ,k_{w};n_{1},n_{2},\cdots ,n_{w}}\left(
x_{1},x_{2},\cdots ,x_{w};q\right) &=&\dprod\limits_{i=1}^{w}\left[ \binom{%
n_{i}-1}{k_{i}}+\binom{n_{i}-1}{k_{i}-1}\right]
[x_{i}]_{q}^{k_{i}}[1-x_{i}]_{q}^{n_{i}-k_{i}} \\
&=&\dprod\limits_{i=1}^{w}\left( [1-x_{i}]_{q}B_{k_{i};n_{i}-1}\left(
x_{i};q\right) +[x_{i}]_{q}B_{k_{i}-1;n_{i}-1}\left( x_{i};q\right) \right) 
\text{.}
\end{eqnarray*}%
This is the desired result.
\end{proof}
\end{theorem}

\begin{remark}
By setting $w=1$ and $q\rightarrow 1^{-}$ into (\ref{Equation 60}), we get
the familiar identity for $B_{k_{1},n_{1}}\left( x_{1}\right) $ as follows: 
\begin{equation*}
B_{k_{1},n_{1}}\left( x_{1}\right) =\left( 1-x_{1}\right)
B_{k_{1},n_{1}-1}\left( x_{1}\right) +x_{1}B_{k_{1}-1,n_{1}-1}\left(
x_{1}\right) \text{.}
\end{equation*}%
(see \cite{Acikgoz1},\cite{Acikgoz3},\cite{SimsekAcikgoz1}).
\end{remark}

\begin{theorem}
For $k_{i},n_{i}\in 
\mathbb{N}
_{0}$ and $x_{i}\in \mathcal{D}$ with $i=1,2,...,w$, we have 
\begin{gather}
B_{n_{1}-k_{1},n_{2}-k_{2}\cdots ,n_{w}-k_{w};n_{1},n_{2},\cdots
,n_{w}}\left( 1-x_{1},1-x_{2},\cdots ,1-x_{w};q\right)  \notag \\
=B_{k_{1},k_{2},\cdots ,k_{w};n_{1},n_{2},\cdots ,n_{w}}\left(
x_{1},x_{2},\cdots ,x_{w};q\right)  \label{Equation151}
\end{gather}
\end{theorem}

\begin{remark}
By substituting $w=1$ and $q\rightarrow 1^{-}$ into (\ref{Equation151}), we
get the well-known identity as follows:%
\begin{equation*}
B_{n_{1}-k_{1},n_{1}}\left( 1-x_{1}\right) =B_{k_{1},n_{1}}\left(
x_{1}\right) \text{.}
\end{equation*}%
(see \cite{Acikgoz1},\cite{Acikgoz3}).
\end{remark}

\begin{definition}
Let $f$ be a continuous function of several variables on $\mathcal{D}^{w}$.
Then the modified $q$-Bernstein operator of order $\dsum%
\limits_{i=1}^{w}n_{i}$ for $f$ is defined by 
\begin{gather}
\mathcal{B}_{n_{1},n_{2},\cdots ,n_{w}}\left( f:x_{1},x_{2},\cdots
,x_{w};q\right)  \label{Equation160} \\
=\dsum\limits_{k_{1}=0}^{n_{1}}\dsum\limits_{k_{2}=0}^{n_{2}}\cdots
\dsum\limits_{k_{W}=0}^{n_{W}}f\left( \frac{k_{1}}{n_{1}},\frac{k_{2}}{n_{2}}%
,\cdots ,\frac{k_{w}}{n_{w}}\right) B_{k_{1},k_{2},\cdots
,k_{w};n_{1},n_{2},\cdots ,n_{w}}\left( x_{1},x_{2},\cdots ,x_{w};q\right) 
\notag
\end{gather}%
where $x_{i}\in \mathcal{D}$, $n_{i}\in 
\mathbb{N}
.$
\end{definition}

\bigskip When we set $f\left( \frac{k_{1}}{n_{1}},\frac{k_{2}}{n_{2}},\cdots
,\frac{k_{w}}{n_{w}}\right) =1$ into (\ref{Equation160}), we easily see
that, 
\begin{gather}
\mathcal{B}_{n_{1},n_{2},\cdots ,n_{w}}\left( 1:x_{1},x_{2},\cdots
,x_{w};q\right)  \label{Equation 100} \\
=\dsum\limits_{k_{1}=0}^{n_{1}}\dsum\limits_{k_{2}=0}^{n_{2}}\cdots
\dsum\limits_{k_{w}=0}^{n_{w}}B_{k_{1},k_{2},\cdots
,k_{w};n_{1},n_{2},\cdots ,n_{w}}\left( x_{1},x_{2},\cdots ,x_{w};q\right) 
\notag
\end{gather}

From the definition of binomial theorem and (\ref{Equation 100}), we get the
following Corollary \ref{Corollary 2.8.} for modified $q$-Bernstein
polynomials for functions of several variables:

\begin{corollary}
\label{Corollary 2.8.}For any $k_{i},n_{i}\in 
\mathbb{N}
_{0}$ and $x_{i}\in \mathcal{D}$ with $i=1,2,\cdots ,w$, we have 
\begin{equation}
\mathcal{B}_{n_{1},n_{2},\cdots ,n_{w}}\left( 1:x_{1},x_{2},\cdots
,x_{w};q\right) =\dprod\limits_{i=1}^{w}\left( 1+\left( 1-q\right)
[x_{i}]_{q}[1-x_{i}]_{q}\right) ^{n_{i}}\text{,}  \label{Equation199}
\end{equation}%
we easily see that 
\begin{equation*}
\underset{q\rightarrow 1}{\lim }\mathcal{B}_{n_{1},n_{2},\cdots
,n_{w}}\left( 1:x_{1},x_{2},\cdots ,x_{w};q\right) =1\text{.}
\end{equation*}%
This is a partition of unity for modified Bernstein polynomials for
functions of several variables.
\end{corollary}

\begin{theorem}
For $\xi _{j}\in 
\mathbb{C}
$, $x_{j}\in \mathcal{D}$ and $n_{j}\in 
\mathbb{N}
$, with $j=1,2,\cdots ,w$ and $i=\sqrt{-1},$ we have 
\begin{equation}
B_{k_{1},k_{2},\cdots ,k_{w};n_{1},n_{2},\cdots ,n_{w}}\left(
x_{1},x_{2},\cdots ,x_{w};q\right) =\frac{1}{\left( 2\pi i\right) ^{w}}%
\underset{w\text{-times}}{\underbrace{\doint\limits_{C}\doint\limits_{C}%
\cdots \doint\limits_{C}}}\dprod\limits_{j=1}^{w}n_{j}!F_{q}^{\left(
k_{j}\right) }\left( x_{j},\xi _{j}\right) \frac{d\xi _{j}}{\xi
_{j}^{n_{j}+1}}  \label{Equa1400}
\end{equation}%
where 
\begin{equation*}
F_{q}^{\left( k\right) }\left( x,t\right) =\frac{\left( t[x]_{q}\right) ^{k}%
}{k!}e^{t[1-x]_{q}}\text{ (see \cite{SimsekAcikgoz1}) }
\end{equation*}%
and $C$ is a circle around the origin and integration is in the positive
direction.

\begin{proof}
By using the definition of the modified $q$-Bernstein polynomials of several
variables and the basic theory of complex analysis including Laurent series
that%
\begin{eqnarray}
&&\underset{w\text{-times}}{\underbrace{\doint\limits_{C}\doint\limits_{C}%
\cdots \doint\limits_{C}}}\dprod\limits_{j=1}^{w}F_{q}^{\left( k_{j}\right)
}\left( x_{j},\xi _{j}\right) \frac{d\xi _{j}}{\xi _{j}^{n_{j}+1}}  \notag \\
&=&\dsum\limits_{l_{1}=0}^{\infty }\dsum\limits_{l_{2}=0}^{\infty }\cdots
\dsum\limits_{l_{w}=0}^{\infty }\doint\limits_{C}\doint\limits_{C}\cdots
\doint\limits_{C}\dprod\limits_{j=1}^{w}\frac{B_{k_{j},l_{j}}\left(
x_{j},q\right) \xi _{j}^{l_{j}}}{l_{j}!}\frac{d\xi _{j}}{\xi _{j}^{n_{j}+1}}
\notag \\
&=&\left( 2\pi i\right) ^{w}\left( \frac{B_{k_{1},k_{2},\cdots
,k_{w};n_{1},n_{2},\cdots ,n_{w}}\left( x_{1},x_{2},\cdots ,x_{w};q\right) }{%
n_{1}!n_{2}!\cdots n_{w}!}\right) .  \label{Equation152}
\end{eqnarray}%
By using (\ref{Equation152}), we obtain 
\begin{eqnarray*}
&&B_{k_{1},k_{2},\cdots ,k_{w};n_{1},n_{2},\cdots ,n_{w}}\left(
x_{1},x_{2},\cdots ,x_{w};q\right) \\
&=&\frac{1}{\left( 2\pi i\right) ^{w}}\underset{w\text{-times}}{\underbrace{%
\doint\limits_{C}\doint\limits_{C}\cdots \doint\limits_{C}}}%
\dprod\limits_{j=1}^{w}n_{j}!F_{q}^{\left( k_{j}\right) }\left( x_{j},\xi
_{j}\right) \frac{d\xi _{j}}{\xi _{j}^{n_{j}+1}}
\end{eqnarray*}%
and 
\begin{equation}
\underset{w\text{-times}}{\underbrace{\doint\limits_{C}\doint\limits_{C}%
\cdots \doint\limits_{C}}}\dprod\limits_{j=1}^{w}F_{q}^{\left( k_{j}\right)
}\left( x_{j},\xi _{j}\right) \frac{d\xi _{j}}{\xi _{j}^{n_{j}+1}}=\left(
2\pi i\right) ^{w}\left( \dprod\limits_{j=1}^{w}\frac{\left[ x_{j}\right]
_{q}^{k_{j}}\left[ 1-x_{j}\right] _{q}^{n_{j}-k_{j}}}{k_{j}!\left(
n_{j}-k_{j}\right) !}\right) .  \label{Equation154}
\end{equation}%
We also obtain from (\ref{Equation151}) and (\ref{Equation154}) that 
\begin{equation}
\frac{1}{\left( 2\pi i\right) ^{w}}\underset{w\text{-times}}{\underbrace{%
\doint\limits_{C}\doint\limits_{C}\cdots \doint\limits_{C}}}%
\dprod\limits_{j=1}^{w}n_{j}!F_{q}^{\left( k_{j}\right) }\left( x_{j},\xi
_{j}\right) \frac{d\xi _{j}}{\xi _{j}^{n_{j}+1}}=\dprod\limits_{j=1}^{w}%
\binom{n_{j}}{k_{j}}[x_{j}]_{q}^{k_{j}}[1-x_{j}]_{q}^{n_{j}-k_{j}}.
\label{Equation155}
\end{equation}%
So, from (\ref{Equation152}) and (\ref{Equation155}) and Corollary \ref%
{Corollary 2.2.}, we complete the proof of theorem.
\end{proof}
\end{theorem}

We now give the modified $q$-Bernstein polynomials for functions of several
variables as a linear combination of polynomials of higher order as follows:

\begin{theorem}
For $k_{i},n_{i}\in 
\mathbb{N}
_{0},$ $x_{i}\in \mathcal{D}$, and $i=1,2,...,w$, we have%
\begin{gather*}
B_{k_{1},k_{2},\cdots ,k_{w};n_{1},n_{2},\cdots ,n_{w}}\left(
x_{1},x_{2},\cdots ,x_{w};q\right) \\
=\dprod\limits_{i=1}^{w}\left[ \left( \frac{n_{i}-k_{i}+1}{k_{i}}\right) 
\frac{[x_{i}]_{q}}{[1-x_{i}]_{q}}\right] B_{k_{1}-1,k_{2}-1,\cdots
,k_{w}-1;n_{1},n_{2},\cdots ,n_{w}}\left( x_{1},x_{2},\cdots ,x_{w};q\right)
.
\end{gather*}

\begin{proof}
Using the definition of modified $q$-Bernstein polynomials for functions of
several variables and the property (\ref{Equation59}), the proof follows.
\end{proof}
\end{theorem}

\begin{theorem}
If $n_{i},k_{i}\in 
\mathbb{N}
_{0}$ and $x_{i}\in \mathcal{D}$ with $i=1,2,...,w$, we have%
\begin{gather*}
B_{k_{1},k_{2},\cdots ,k_{w};n_{1},n_{2},\cdots ,n_{w}}\left(
x_{1},x_{2},\cdots ,x_{w};q\right) \\
=\dsum\limits_{l_{1}=k_{1}}^{n_{1}}\dsum\limits_{l_{2}=k_{2}}^{n_{2}}\cdots
\dsum\limits_{l_{w}=k_{w}}^{n_{w}}\dprod\limits_{i=1}^{w}\binom{n_{i}}{l_{i}}%
\binom{l_{i}}{k_{i}}\left( -1\right) ^{l_{i}-k_{i}}q^{\left(
l_{i}-k_{i}\right) \left( 1-x_{i}\right) }[x_{i}]_{q}^{l_{i}}.
\end{gather*}

\begin{proof}
From the definition of modified $q$-Bernstein polynomials of several
variables and binomial theorem with $n_{i},$ $k_{i}\in 
\mathbb{N}
_{0}$ and $x_{i}\in \mathcal{D}$\ for $i=1,2,...,w,$ we have%
\begin{gather*}
B_{k_{1},k_{2},\cdots ,k_{w};n_{1},n_{2},\cdots ,n_{w}}\left(
x_{1},x_{2},\cdots ,x_{w};q\right) =\dprod\limits_{i=1}^{w}\binom{n_{i}}{%
k_{i}}[x_{i}]_{q}^{k_{i}}[1-x_{i}]_{q}^{n_{i}-k_{i}} \\
=\dsum\limits_{l_{1}=k_{1}}^{n_{1}}\dsum\limits_{l_{2}=k_{2}}^{n_{2}}\cdots
\dsum\limits_{l_{w}=k_{w}}^{n_{w}}\dprod\limits_{i=1}^{w}\binom{n_{i}}{l_{i}}%
\binom{l_{i}}{k_{i}}\left( -1\right) ^{l_{i}-k_{i}}q^{\left(
l_{i}-k_{i}\right) \left( 1-x_{i}\right) }[x_{i}]_{q}^{l_{i}}.
\end{gather*}%
This is the desired result.
\end{proof}
\end{theorem}

\begin{theorem}
For $n_{i},l_{i}\in 
\mathbb{N}
_{0}$ and $x_{i}\in \mathcal{D}$, with $i=1,2,...,w,$ we have%
\begin{gather*}
\left( \dprod\limits_{i=1}^{w}[x_{i}]_{q}\right) ^{m} \\
=\dprod\limits_{i=1}^{w}\frac{1}{\left( \left[ 1-x_{i}\right] _{q}+\left[
x_{i}\right] _{q}\right) ^{n_{i}-m}}\dsum\limits_{k_{1}=m}^{n_{1}}\dsum%
\limits_{k_{2}=m}^{n_{2}}\cdots
\dsum\limits_{k_{w}=m}^{n_{w}}\dprod\limits_{i=1}^{w}\frac{\binom{k_{i}}{m}}{%
\binom{n_{i}}{m}}B_{k_{1},k_{2},\cdots ,k_{w};n_{1},n_{2},\cdots
,n_{w}}\left( x_{1},x_{2},\cdots ,x_{w};q\right) .
\end{gather*}

\begin{proof}
We easily see from the property of the modified $q$-Bernstein polynomials of
several variables that%
\begin{gather*}
\dsum\limits_{k_{1}=1}^{n_{1}}\dsum\limits_{k_{2}=1}^{n_{2}}\cdots
\dsum\limits_{k_{w}=1}^{n_{w}}\dprod\limits_{i=1}^{w}\frac{k_{i}}{n_{i}}%
B_{k_{1},k_{2},\cdots ,k_{w};n_{1},n_{2},\cdots ,n_{w}}\left(
x_{1},x_{2},\cdots ,x_{w};q\right) \\
=\dprod\limits_{i=1}^{w}\left[ x_{i}\right] _{q}\left( \left[ x_{i}\right]
_{q}+\left[ 1-x_{i}\right] _{q}\right) ^{n_{i}-1}
\end{gather*}%
and also 
\begin{gather*}
\dsum\limits_{k_{1}=2}^{n_{1}}\dsum\limits_{k_{2}=2}^{n_{2}}\cdots
\dsum\limits_{k_{w}=2}^{n_{w}}\dprod\limits_{i=1}^{w}\frac{\binom{k_{i}}{2}}{%
\binom{n_{i}}{2}}B_{k_{1},k_{2},\cdots ,k_{w};n_{1},n_{2},\cdots
,n_{w}}\left( x_{1},x_{2},\cdots ,x_{w};q\right) \\
=\left( \dprod\limits_{i=1}^{w}\left[ x_{i}\right] _{q}\right) ^{2}\left( %
\left[ x_{i}\right] _{q}+\left[ 1-x_{i}\right] _{q}\right) ^{n_{i}-2}\text{.}
\end{gather*}%
Continuing this method, we have 
\begin{gather*}
\left( \dprod\limits_{i=1}^{w}[x_{i}]_{q}\right) ^{m}=\dprod\limits_{i=1}^{w}%
\frac{1}{\left( \left[ 1-x_{i}\right] _{q}+\left[ x_{i}\right] _{q}\right)
^{n_{i}-m}} \\
\times \dsum\limits_{k_{1}=m}^{n_{1}}\dsum\limits_{k_{2}=m}^{n_{2}}\cdots
\dsum\limits_{k_{w}=m}^{n_{w}}\dprod\limits_{i=1}^{w}\frac{\binom{k_{i}}{m}}{%
\binom{n_{i}}{m}}B_{k_{1},k_{2},\cdots ,k_{w};n_{1},n_{2},\cdots
,n_{w}}\left( x_{1},x_{2},\cdots ,x_{w};q\right)
\end{gather*}%
and after making some algebraic operations, we obtain the desired result.
\end{proof}
\end{theorem}

We have seen from the theorem given above, it is possible to write $\left(
\dprod\limits_{i=1}^{w}[x_{i}]_{q}\right) ^{m}$ as a linear combination of
modified $q$-Bernstein polynomials of several variables by using the degree
evaluation formulae and mathematical induction method.

For $k\in 
\mathbb{N}
_{0}$, the Bernoulli polynomials of degree $k$ are defined by 
\begin{eqnarray*}
\underset{k-times}{\underbrace{\left( \frac{t}{e^{t}-1}\right) \left( \frac{t%
}{e^{t}-1}\right) \times \cdots \times \left( \frac{t}{e^{t}-1}\right) }}%
e^{xt} &=&\left( \frac{t}{e^{t}-1}\right) ^{k}e^{xt} \\
&=&\dsum\limits_{n=0}^{\infty }B_{n}^{\left( k\right) }\left( x\right) \frac{%
t^{n}}{n!},
\end{eqnarray*}%
and $B_{n}^{\left( k\right) }=B_{n}^{\left( k\right) }\left( 0\right) $ are
called the $n$-th Bernoulli numbers of order $k$. It is well known that the
second kind Stirling numbers are defined by $\frac{\left( e^{t}-1\right) ^{k}%
}{k!}:=\dsum\limits_{n=0}^{\infty }S\left( n,k\right) \frac{t^{n}}{n!}$ for $%
k\in 
\mathbb{N}
$ (see \cite{Kim-Lee-Chae Jang},\cite{SimsekAcikgoz1}). By using the above
relations, we can give the following theorem:

\begin{theorem}
For $k_{i},n_{i},\in 
\mathbb{N}
_{0}$ and $x_{i}\in \mathcal{D}$ with $i=1,2,\cdots ,w$, we have%
\begin{gather*}
B_{k_{1},k_{2},\cdots ,k_{w};n_{1},n_{2},\cdots ,n_{w}}\left(
x_{1},x_{2},\cdots ,x_{w};q\right) \\
=\dsum\limits_{l_{1}=0}^{n_{1}}\dsum\limits_{l_{2}=0}^{n_{2}}\cdots
\dsum\limits_{l_{w}=0}^{n_{w}}\dprod\limits_{i=1}^{w}[x_{i}]_{q}^{l_{i}}%
\binom{n_{i}}{l_{i}}B_{l_{i}}^{\left( k_{i}\right) }\left( \left[ 1-x_{i}%
\right] _{q}\right) S\left( n_{i}-l_{i},k_{i}\right) .
\end{gather*}

\begin{proof}
By using the generating function of modified $q$-Bernstein polynomials of
several variables, we have 
\begin{gather*}
\dprod\limits_{i=1}^{w}\frac{\left( t[x_{i}]_{q}\right) ^{k_{i}}}{k_{i}!}%
e^{t\left( \dsum\limits_{i=1}^{w}[1-x_{i}]_{q}\right)
}=\dprod\limits_{i=1}^{w}[x_{i}]_{q}^{k_{i}}\left(
\dsum\limits_{n_{1}=0}^{\infty }S\left( n_{1},k_{1}\right) \frac{t^{n_{1}}}{%
n_{1}!}\right) \cdots \left( \dsum\limits_{n_{w}=0}^{\infty }S\left(
n_{w},k_{w}\right) \frac{t^{n_{w}}}{n_{w}!}\right) \times \\
\left( \dsum\limits_{l_{1}=0}^{\infty }B_{l_{1}}^{\left( k_{1}\right)
}\left( [1-x_{1}]_{q}\right) \frac{t^{l_{1}}}{l_{1}!}\right) \cdots \left(
\dsum\limits_{l_{w}=0}^{\infty }B_{l_{w}}^{\left( k_{w}\right) }\left(
[1-x_{w}]_{q}\right) \frac{t^{l_{w}}}{l_{w}!}\right) .
\end{gather*}%
By using the Cauchy product for sums given above 
\begin{gather*}
B_{k_{1},k_{2},\cdots ,k_{w};n_{1},n_{2},\cdots ,n_{w}}\left(
x_{1},x_{2},\cdots ,x_{w};q\right) \\
=\dsum\limits_{l_{1}=0}^{n_{1}}\dsum\limits_{l_{2}=0}^{n_{2}}\cdots
\dsum\limits_{l_{w}=0}^{n_{w}}\dprod\limits_{i=1}^{w}[x_{i}]_{q}^{l_{i}}%
\binom{n_{i}}{l_{i}}B_{l_{i}}^{\left( k_{i}\right) }\left( \left[ 1-x_{i}%
\right] _{q}\right) S\left( n_{i}-l_{i},k_{i}\right) .
\end{gather*}%
By comparing the last two relations, we have the desired result.
\end{proof}
\end{theorem}

Let $\Delta $ be the shift difference operator defined by $\Delta f\left(
x\right) =f\left( x+1\right) -f\left( x\right) $. By using the mathematical
induction method we have%
\begin{equation}
\Delta ^{n}f\left( 0\right) =\dsum\limits_{k=0}^{n}\binom{n}{k}\left(
-1\right) ^{n-k}f\left( k\right) ,  \label{Equation90}
\end{equation}%
for $n\in 
\mathbb{N}
$ and using (\ref{Equation90}) in the generating function of second kind
Stirling numbers,%
\begin{gather}
\dsum\limits_{n=0}^{\infty }S\left( n,k\right) \frac{t^{n}}{n!}=\frac{1}{k!}%
\dsum\limits_{l=0}^{k}\binom{k}{l}\left( -1\right) ^{k-l}e^{lt}  \notag \\
=\dsum\limits_{n=0}^{\infty }\left( \frac{1}{k!}\dsum\limits_{l=0}^{k}\binom{%
k}{l}\left( -1\right) ^{k-l}l^{n}\right) \frac{t^{n}}{n!}\text{.}
\label{Equation150}
\end{gather}

By comparing the coefficients on both sides, we have 
\begin{equation}
S\left( n,k\right) =\frac{1}{k!}\dsum\limits_{l=0}^{k}\binom{k}{l}\left(
-1\right) ^{k-l}l^{n}\text{.}  \label{Equation91}
\end{equation}

When we compared Eq. (\ref{Equation90}) and Eq. (\ref{Equation91}), becomes%
\begin{equation}
S\left( n,k\right) =\frac{\Delta ^{k}0^{n}}{k!}\text{.}  \label{Equation92}
\end{equation}

For $n_{i},k_{i}\in 
\mathbb{N}
$, by using the equation (\ref{Equation92}), we obtain the relation 
\begin{gather*}
B_{k_{1},k_{2},\cdots ,k_{w};n_{1},n_{2},\cdots ,n_{w}}\left(
x_{1},x_{2},\cdots ,x_{w};q\right) \\
=\dsum\limits_{l_{1}=0}^{n_{1}}\dsum\limits_{l_{2}=0}^{n_{2}}\cdots
\dsum\limits_{l_{w}=0}^{n_{w}}\dprod\limits_{i=1}^{w}[x_{i}]_{q}^{l_{i}}%
\binom{n_{i}}{l_{i}}B_{l_{i}}^{\left( k_{i}\right) }\left( \left[ 1-x_{i}%
\right] _{q}\right) \frac{\Delta ^{k_{i}}0^{n_{i}-l_{i}}}{k_{i}!}
\end{gather*}%
which is the relation of the $q$-Bernstein polynomials of several variables
in terms of Bernoulli polynomials of order $k$ and second Stirling numbers
with shift difference operator.

Let $\left( Eh\right) \left( x\right) =h\left( x+1\right) $ be the shift
operator. Then the $q$-difference operator is defined by 
\begin{equation}
\Delta _{q}^{n}=\dprod\limits_{i=0}^{n-1}\left( E-q^{i}I\right)
\label{Equation93}
\end{equation}%
where $I$ is the identity operator (see \cite{Kim-Lee-Chae Jang}).

For $f\in C\left( [0,1]\right) $ and $n\in 
\mathbb{N}
$, we have 
\begin{equation}
\Delta _{q}^{n}f\left( 0\right) =\dsum\limits_{k=0}^{n}\binom{n}{k}%
_{q}\left( -1\right) ^{k}q^{\binom{n}{2}}f\left( n-k\right) \text{,}
\label{Equation94}
\end{equation}%
where $\binom{n}{k}_{q}$ is the Gaussian binomial coefficient defined by 
\begin{equation}
\binom{n}{k}_{q}=\frac{[n]_{q}[n-1]_{q}\cdots \lbrack n-k+1]_{q}}{[k]_{q}!}.
\label{Equation95}
\end{equation}

\begin{theorem}
For $n_{i},l_{i}\in 
\mathbb{N}
_{0}$ and $x_{i}\in \mathcal{D}$ for $i=1,2,...,w$, we have%
\begin{gather*}
\dprod\limits_{i=1}^{w}\frac{1}{\left( \left[ 1-x_{i}\right] _{q}+\left[
x_{i}\right] _{q}\right) ^{n_{i}-l_{i}}}\dsum\limits_{k_{1}=m}^{n_{1}}\dsum%
\limits_{k_{2}=m}^{n_{2}}\cdots \dsum\limits_{k_{w}=m}^{n_{w}}\left(
\dprod\limits_{i=1}^{w}\frac{\binom{k_{i}}{m}}{\binom{n_{i}}{m}}\right) \\
\times B_{k_{1},k_{2},\cdots ,k_{w};n_{1},n_{2},\cdots ,n_{w}}\left(
x_{1},x_{2},\cdots ,x_{w};q\right) \\
=\dsum\limits_{l_{1}=0}^{m}\dsum\limits_{l_{2}=0}^{m}\cdots
\dsum\limits_{l_{w}=0}^{m}q^{\dsum\limits_{i=1}^{w}\binom{l_{i}}{2}%
}\dprod\limits_{i=1}^{w}\binom{x_{i}}{l_{i}}\left[ l_{i}\right] _{q}!S\left(
m,l_{i};q\right) \text{.}
\end{gather*}

\begin{proof}
To prove this theorem, we let $F_{q}\left( t\right) $ be the generating
function of the $q$-extension of the second kind Stirling numbers as follows:%
\begin{equation}
F_{q}\left( t\right) :=\frac{q^{-\binom{k}{2}}}{[k]_{q}!}\dsum%
\limits_{i=0}^{k}\left( -1\right) ^{k-i}\binom{k}{i}_{q}q^{\binom{k-i}{2}%
}e^{[i]_{q}t}=\dsum\limits_{n=0}^{\infty }S\left( n,k;q\right) \frac{t^{n}}{%
n!}  \notag
\end{equation}%
From the above, we have 
\begin{eqnarray}
S\left( n,k;q\right) &=&\frac{q^{-\binom{k}{2}}}{[k]_{q}!}%
\dsum\limits_{i=0}^{k}\left( -1\right) ^{i}q^{\binom{i}{2}}\binom{k}{i}%
_{q}[k-i]_{q}^{n}  \notag \\
&=&\frac{q^{-\binom{k}{2}}}{[k]_{q}!}\Delta _{q}^{k}0^{n}  \label{Equation97}
\end{eqnarray}%
where $[k]_{q}!=[k]_{q}[k-1]_{q}\cdots \lbrack 2]_{q}[1]_{q}.$ It is easy to
see that 
\begin{equation}
\lbrack x]_{q}^{n}=\dsum\limits_{k=0}^{n}q^{\binom{k}{2}}\binom{x}{k}%
_{q}[k]_{q}!S\left( n,k;q\right)  \label{Equation98}
\end{equation}%
and in similar way that 
\begin{equation}
\left( \dprod\limits_{i=1}^{w}[x_{i}]_{q}\right)
^{m}=\dsum\limits_{l_{1}=0}^{m}\dsum\limits_{l_{2}=0}^{m}\cdots
\dsum\limits_{l_{w}=0}^{m}q^{\dsum\limits_{i=1}^{w}\binom{l_{i}}{2}%
}\dprod\limits_{i=1}^{w}\binom{x_{i}}{l_{i}}\left[ l_{i}\right] _{q}!S\left(
m,l_{i};q\right) \text{.}  \label{Equation99}
\end{equation}%
Then, we obtain the desired result from (\ref{Equation98}) and (\ref%
{Equation99}).
\end{proof}
\end{theorem}

\section{\textbf{Interpolation Function of Modified }$q$\textbf{-Bernstein
Polynomials for Functions of Several Variables}}

The classical Bernoulli numbers interpolate by Riemann zeta function, which
has profound effect on Analytic numbers theory and complex analysis. The
values of the negative integer points, also found by Euler, are rational
numbers and play a vital and important role in the theory of modular forms.
Many generalization of the Riemann zeta function, such as Dirichlet series,
Dirichlet $L$-functions and $L$-functions, are known in \cite{Kim 5}, \cite%
{Choi 1}, \cite{Choi 2}, \cite{Choi 3}, \cite{Araci2}, \cite{Araci3}. So, we
construct interpolation function of modified $q$-Bernstein polynomials of
several variables.

For $s\in 
\mathbb{C}
$ and $x_{i}\neq 1$ with $i=1,2,...,w$, by applying Mellin transformation to
Eq. (\ref{Equation58}), we procure%
\begin{align*}
& D_{q}\left( s,k_{1},k_{2}\cdots k_{w};x_{1},x_{2}\cdots x_{w}\right) \\
& =\frac{1}{\Gamma \left( s\right) }\int_{0}^{\infty }t^{s-\left(
k_{1}+k_{2}+\cdots +k_{w}\right) -1}F_{k_{1},k_{2},\cdots ,k_{w}}\left(
-t,q;x_{1},x_{2},\cdots ,x_{w}\right) dt \\
& =\frac{1}{\Gamma \left( s\right) }\int_{0}^{\infty }t^{s-\left(
k_{1}+k_{2}+\cdots +k_{w}\right) -1}\dprod\limits_{i=1}^{w}\frac{\left( -t%
\left[ x_{i}\right] _{q}\right) ^{k_{i}}}{k_{i}!}e^{-t\dsum\limits_{i=1}^{w}%
\left[ 1-x_{i}\right] _{q}}dt \\
& =\left( -1\right) ^{\dsum\limits_{i=1}^{w}k_{i}}\dprod\limits_{i=1}^{w}%
\frac{\left[ x_{i}\right] _{q}^{k_{i}}}{k_{i}!}\left( \frac{1}{\Gamma \left(
s\right) }\int_{0}^{\infty }t^{s-1}e^{-t\left[ 1-x_{i}\right] _{q}}dt\right)
\\
& =\left( -1\right) ^{\dsum\limits_{i=1}^{w}k_{i}}\dprod\limits_{i=1}^{w}%
\frac{\left[ x_{i}\right] _{q}^{k_{i}}}{k_{i}!}\left[ 1-x_{i}\right]
_{q}^{-s}.
\end{align*}

From the above, we give the definition of interpolation function for
Corollary \ref{Corollary 2.2.} as follows:

\begin{definition}
Let $s\in 
\mathbb{C}
$ and $x_{i}\neq 1$ with $i=1,2,...,w$. We define interpolation function of
the polynomials $B_{k_{1},k_{2},\cdots ,k_{w};n_{1},n_{2},\cdots
,n_{w}}\left( x_{1},x_{2},\cdots ,x_{w};q\right) $ as 
\begin{equation}
D_{q}\left( s,k_{1},k_{2}\cdots k_{w};x_{1},x_{2}\cdots x_{w}\right) =\left(
-1\right) ^{\dsum\limits_{i=1}^{w}k_{i}}\dprod\limits_{i=1}^{w}\frac{\left[
x_{i}\right] _{q}^{k_{i}}}{k_{i}!}\left( \left[ 1-x_{i}\right] _{q}\right)
^{-s}\text{.}  \label{Equation50}
\end{equation}
\end{definition}

\begin{remark}
By substituting $w=1$ into (\ref{Equation50}), we get%
\begin{equation*}
D_{q}\left( s,k_{1}\right) =\left( -1\right) ^{k_{1}}\frac{\left[ x_{1}%
\right] _{q}^{k_{1}}}{k_{1}!}\left[ 1-x_{1}\right] _{q}^{-s}
\end{equation*}%
where $D_{q}\left( s,k_{1}\right) $ is introduced by Simsek and Acikgoz cf. 
\cite{SimsekAcikgoz1}.
\end{remark}

\bigskip Substituting $s=-\left( n_{1}+n_{2}+\cdots +n_{w}\right) $ into Eq.
(\ref{Equation50}), we have%
\begin{align*}
& D_{q}\left( -n_{1}-n_{2}-\cdots -n_{w},k_{1},k_{2}\cdots
k_{w};x_{1},x_{2}\cdots x_{w}\right) \\
& =\left( -1\right) ^{\dsum\limits_{i=1}^{w}k_{i}}\dprod\limits_{i=1}^{w}%
\frac{\left[ x_{i}\right] _{q}^{k_{i}}}{k_{i}!}\left[ 1-x_{i}\right]
_{q}^{n_{i}} \\
& =\dprod\limits_{i=1}^{w}\frac{\left( -1\right) ^{k_{i}}n_{i}!}{\left(
n_{i}+k_{i}\right) !}\dprod\limits_{i=1}^{w}\binom{n_{i}+k_{i}}{k_{i}}\left[
x_{i}\right] _{q}^{k_{i}}\left[ 1-x_{i}\right] _{q}^{\left(
n_{i}+k_{i}\right) -k_{i}} \\
& =\dprod\limits_{i=1}^{w}\frac{\left( -1\right) ^{k_{i}}n_{i}!}{\left(
n_{i}+k_{i}\right) !}B_{k_{1},k_{2},\cdots
,k_{w};n_{1}+k_{1},n_{2}+k_{2},\cdots ,n_{w}+k_{w}}\left( x_{1},x_{2},\cdots
,x_{w};q\right) \text{.}
\end{align*}

So, we arrive at the following theorem.

\begin{theorem}
The following equality holds true:%
\begin{gather*}
D_{q}\left( -n_{1}-n_{2}-\cdots -n_{w},k_{1},k_{2}\cdots
k_{w};x_{1},x_{2}\cdots x_{w}\right) \\
=\dprod\limits_{i=1}^{w}\frac{\left( -1\right) ^{k_{i}}n_{i}!}{\left(
n_{i}+k_{i}\right) !}B_{k_{1},k_{2},\cdots
,k_{w};n_{1}+k_{1},n_{2}+k_{2},\cdots ,n_{w}+k_{w}}\left( x_{1},x_{2},\cdots
,x_{w};q\right) \text{.}
\end{gather*}
\end{theorem}

By using (\ref{Equation50}), we have 
\begin{equation*}
D_{q}\left( s,k_{1},k_{2}\cdots k_{w};x_{1},x_{2}\cdots x_{w}\right)
\rightarrow D\left( s,k_{1},k_{2}\cdots k_{w};x_{1},x_{2}\cdots x_{w}\right) 
\text{ as }q\rightarrow 1.
\end{equation*}

Thus one has%
\begin{equation}
D\left( s,k_{1},k_{2}\cdots k_{w};x_{1},x_{2}\cdots x_{w}\right) =\left(
-1\right) ^{\dsum\limits_{i=1}^{w}k_{i}}\dprod\limits_{i=1}^{w}\frac{%
x_{i}^{k_{i}}}{k_{i}!}\left( 1-x_{i}\right) ^{-s}\text{.}  \label{Equation51}
\end{equation}

By substituting $x_{i}=1$ with $i=1,2,...,w$ within the above, we have 
\begin{equation*}
D\left( s,k_{1},k_{2}\cdots k_{w};x_{1},x_{2}\cdots x_{w}\right) =\infty .
\end{equation*}

We now evaluate the $i$-th $s$-derivative of $D\left( s,k_{1},k_{2}\cdots
k_{w};x_{1},x_{2}\cdots x_{w}\right) $ as follows: For $x_{j}\neq 1$ with $%
i=1,2,...,w$%
\begin{equation}
\frac{\partial ^{i}}{\partial s^{i}}D\left( s,k_{1},k_{2}\cdots
k_{w};x_{1},x_{2}\cdots x_{w}\right) =\log ^{i}\left( \frac{1}{1-x_{i}}%
\right) D\left( s,k_{1},k_{2}\cdots k_{w};x_{1},x_{2}\cdots x_{w}\right)
\label{Equation52}
\end{equation}%
which seems to be interesting.

\begin{remark}
By taking $w=1,$ $q\rightarrow 1^{-}$ into (\ref{Equation51}), we arrive at
the following relation which was proved by Simsek and Acikgoz \cite%
{SimsekAcikgoz1},%
\begin{equation*}
\frac{\partial ^{i}}{\partial s^{i}}D\left( s,k_{1};x_{1}\right) =\log
^{i}\left( \frac{1}{1-x_{1}}\right) D\left( s,k_{1};x_{1}\right) .
\end{equation*}
\end{remark}

\section{$p$\textbf{-adic Integral Representation of }$q$\textbf{%
-Bernstein-type polynomials}}

Throughout this section, we will use the following notations: $%
\mathbb{Z}
_{p}$ denotes the ring of $p$-adic rational integers, $%
\mathbb{Q}
_{p}$ denotes the field of $p$-adic rational numbers, $%
\mathbb{C}
_{p}$ denotes the completion of algebraic closure of $%
\mathbb{Q}
_{p}$. Let $v_{p}$ be the normalized exponential valuation of $%
\mathbb{C}
_{p}$ with $\left\vert p\right\vert _{p}=p^{-v_{p}\left( p\right) }=p^{-1}$.
When we mention about $q$-extension, we say that $q$ is considered in many
ways such as an indeterminate, a complex number $q\in 
\mathbb{C}
$, or $p$-adic number $q\in 
\mathbb{C}
_{p}$. If $q\in 
\mathbb{C}
$ we assume that $\left\vert q\right\vert <1$. If $q\in 
\mathbb{C}
_{p}$we normally assume that $\left\vert q-1\right\vert _{p}<p^{-\frac{1}{p-1%
}}$ so that $q^{x}=\exp \left( x\log q\right) $ for $\left\vert x\right\vert
_{p}\leq 1$ \textit{cf. }\cite{Araci1}, \cite{Araci4}, \cite%
{Araci-Acikgoz-Kilicman}, \cite{Araci}, \cite{Araci2}, \cite{Araci3}, \cite%
{Kim 4}, \cite{Kim7}, \cite{Kim6}. Let $UD\left( 
\mathbb{Z}
_{p}\right) $ be the set of uniformly differentiable function. For $f\in
UD\left( 
\mathbb{Z}
_{p}\right) $, the $p$-adic $q$-integral on $%
\mathbb{Z}
_{p}$ was originally defined by Kim \cite{Kim 4} as follows:%
\begin{eqnarray}
I_{q}\left( f\right) &=&\int_{%
\mathbb{Z}
_{p}}f\left( x\right) d\mu _{q}\left( x\right)  \label{Equation 101} \\
&=&\lim_{n\rightarrow \infty }\sum_{x=0}^{p^{n}-1}f\left( x\right) \mu
_{q}\left( x+p^{n}%
\mathbb{Z}
_{p}\right) =\lim_{n\rightarrow \infty }\frac{1}{\left[ p^{n}\right] _{q}}%
\sum_{x=0}^{p^{n}-1}f\left( x\right) q^{x}\text{.}  \notag
\end{eqnarray}

As $q$ tends to $1^{-}$ in (\ref{Equation 101}), we get known identity ($p$%
-adic Volkenborn Integral) as%
\begin{equation*}
\int_{%
\mathbb{Z}
_{p}}f\left( x\right) d\mu _{1}\left( x\right) =\lim_{n\rightarrow \infty }%
\frac{1}{p^{n}}\sum_{x=0}^{p^{n}-1}f\left( x\right) \text{ (see \cite{Kim7}, 
\cite{Kim6}).}
\end{equation*}

As $I_{-q}\left( f\right) =\lim_{q\rightarrow -q}I_{q}\left( f\right) $
symbolically, which yields, for $p$ an odd prime, to   
\begin{equation}
I_{-q}\left( f\right) =\int_{%
\mathbb{Z}
_{p}}f\left( x\right) d\mu _{-q}\left( x\right) =\lim_{n\rightarrow \infty }%
\frac{1}{\left[ p^{n}\right] _{-q}}\sum_{x=0}^{p^{n}-1}\left( -1\right)
^{x}f\left( x\right) q^{x}  \label{Equation 103}
\end{equation}%
is known as fermionic $p$-adic $q$-invariant integral in the $p$-adic
integer ring. And also, letting $q$ to $1^{-}$ in (\ref{Equation 103}), it
reduces to%
\begin{equation}
\int_{%
\mathbb{Z}
_{p}}f\left( x\right) d\mu _{-1}\left( x\right) =\lim_{n\rightarrow \infty
}\sum_{x=0}^{p^{n}-1}\left( -1\right) ^{x}f\left( x\right) 
\label{Equation 104}
\end{equation}%
(see \cite{Kim4}, \cite{Kim7}, \cite{Kim9}).

\bigskip The Bernoulli numbers was generated by the following generating
function: For $t\in 
\mathbb{C}
$ $\left( \text{with }\left\vert t\right\vert <2\pi \right) $ 
\begin{equation*}
\sum_{n=0}^{\infty }B_{n}\frac{t^{n}}{n!}=\frac{t}{e^{t}-1}\text{ (see \cite%
{Kim7}, \cite{Kim6}, \cite{Kim8}, \cite{Would}).}
\end{equation*}

Next, it was shown that the Bernoulli numbers can be generated by $p$-adic
Volkenborn integral as follows%
\begin{equation*}
B_{n}=\int_{%
\mathbb{Z}
_{p}}x^{n}d\mu _{1}\left( x\right) \text{ for }n\in 
\mathbb{Z}
_{+}:=%
\mathbb{N}
\cup \left\{ 0\right\} \text{, where }%
\mathbb{N}
\text{ is the set of natural numbers.}
\end{equation*}

The following may be defined as a new $q$-extension of Bernoulli numbers%
\begin{equation*}
\beta _{n}\left( q\right) =\int_{%
\mathbb{Z}
_{p}}q^{x}\left[ x\right] _{q}^{n}d\mu _{q}\left( x\right) .
\end{equation*}

Observe that 
\begin{equation*}
\lim_{q\rightarrow 1^{-}}\beta _{n}\left( q\right) =B_{n}.
\end{equation*}

Recall that 
\begin{equation*}
\sum_{n=k}^{\infty }B_{k,n}\left( x;q\right) \frac{t^{n}}{n!}=\frac{\left( t%
\left[ x\right] _{q}\right) ^{k}}{k!}e^{t\left[ 1-x\right] _{q}}
\end{equation*}%
is called $q$-Bernstein-type polynomials. From this, we have%
\begin{equation*}
B_{k,n}\left( x;q\right) =\binom{n}{k}\left[ x\right] _{q}^{k}\left[ 1-x%
\right] _{q}^{n-k}\text{.}
\end{equation*}

Throughout this section, we will assume that $x\in \left( 0,1\right) $. So
we can write%
\begin{equation*}
\left[ x\right] _{q}^{k}=\frac{q^{x}B_{k,n}\left( x;q\right) }{\binom{n}{k}%
\left( 1-\left[ x\right] _{q}\right) ^{n-k}}=\frac{q^{x}}{\binom{n}{k}}%
B_{k,n}\left( x;q\right) \sum_{l=0}^{\infty }\binom{n-k+l-1}{l}\left[ x%
\right] _{q}^{l}\text{.}
\end{equation*}

Further%
\begin{equation}
\frac{\binom{n}{k}}{B_{k,n}\left( x;q\right) }=\sum_{l=0}^{\infty }\binom{%
n-k+l-1}{l}q^{x}\left[ x\right] _{q}^{l-k}\text{.}  \label{Equation 102}
\end{equation}

Applying $p$-adic $q$-integral on $%
\mathbb{Z}
_{p}$ in the both sides of (\ref{Equation 102}), it yields to%
\begin{equation*}
\int_{%
\mathbb{Z}
_{p}}\frac{\binom{n}{k}}{B_{k,n}\left( x;q\right) }d\mu _{q}\left( x\right)
=\sum_{l=k}^{\infty }\binom{n-k+l-1}{l}\beta _{l-k}\left( q\right) \text{.}
\end{equation*}

Therefore we get the following theorem.

\begin{theorem}
\label{thm33}For $k=0,1,2,\cdots ,n$ and $n\in 
\mathbb{Z}
_{+}$, we have%
\begin{equation*}
\int_{%
\mathbb{Z}
_{p}}\binom{n}{k}B_{k,n}^{-1}\left( x;q\right) d\mu _{q}\left( x\right)
=\sum_{l=k}^{\infty }\binom{n-k+l-1}{l}\beta _{l-k}\left( q\right)
\end{equation*}%
where $B_{k,n}^{-1}\left( x;q\right) $ is the inverse of $B_{k,n}\left(
x;q\right) $.
\end{theorem}

As $q$ tends to $1^{-}$ in Theorem \ref{thm33}, we have the following
Corollary.

\begin{corollary}
For $k=0,1,2,\cdots ,n$ and $n\in 
\mathbb{Z}
_{+}$, we have%
\begin{equation*}
\int_{%
\mathbb{Z}
_{p}}\binom{n}{k}B_{k,n}^{-1}\left( x\right) d\mu \left( x\right)
=\sum_{l=k}^{\infty }\binom{n-k+l-1}{l}B_{l-k}
\end{equation*}%
where $B_{k,n}^{-1}\left( x\right) $ is the inverse of $B_{k,n}\left(
x\right) $.
\end{corollary}

The generating function of Euler polynomials has the following series
expansion at $t=0:$%
\begin{equation*}
\sum_{n=0}^{\infty }E_{n}\left( x\right) \frac{t^{n}}{n!}=\frac{2}{e^{t}+1}%
e^{xt}\text{ \ \ \ }\left( \left\vert t\right\vert <\pi \right) .
\end{equation*}

The Euler numbers are defined by $E_{n}\left( 1/2\right) =2^{n}E_{n}$. The
Euler polynomials can be generated through Equation (\ref{Equation 104})%
\begin{equation}
E_{n}\left( x\right) =\int_{%
\mathbb{Z}
_{p}}\left( x+y\right) ^{n}d\mu _{-1}\left( y\right)   \label{Equation 105}
\end{equation}%
(for details, see \cite{Araci5}, \cite{Araci}, \cite{Araci2}, \cite{Kim7}, 
\cite{Kim8}, \cite{Kim9}, \cite{Srivastava}, \cite{Srivastava1}, \cite%
{Srivastava2}).

In \cite{Kim9}, Kim defined the following $q$-Euler numbers%
\begin{equation*}
E_{n,q}=\int_{%
\mathbb{Z}
_{p}}q^{-x}\left[ x\right] _{q}^{n}d\mu _{-q}\left( x\right) .
\end{equation*}

It is clear that 
\begin{equation*}
\lim_{q\rightarrow 1^{-}}E_{n,q}=E_{n}\left( 0\right) \text{.}
\end{equation*}

By (\ref{Equation 103}) and (\ref{Equation 102}), we arrive at the following
theorem.

\begin{theorem}
\label{thm34}For $k=0,1,2,\cdots ,n$ and $n\in 
\mathbb{Z}
_{+}$, we have%
\begin{equation*}
\int_{%
\mathbb{Z}
_{p}}\binom{n}{k}B_{k,n}^{-1}\left( x;q\right) d\mu _{-q}\left( x\right)
=\sum_{l=k}^{\infty }\binom{n-k+l-1}{l}E_{l-k,q}
\end{equation*}%
where $B_{k,n}^{-1}\left( x;q\right) $ is the inverse of $B_{k,n}\left(
x;q\right) $.
\end{theorem}

As $q$ tends to $1^{-}$ in Theorem \ref{thm34}, we have the following
Corollary.

\begin{corollary}
For $k=0,1,2,\cdots ,n$ and $n\in 
\mathbb{Z}
_{+}$, we have%
\begin{equation*}
\int_{%
\mathbb{Z}
_{p}}\binom{n}{k}B_{k,n}^{-1}\left( x\right) d\mu _{-1}\left( x\right)
=\sum_{l=k}^{\infty }\binom{n-k+l-1}{l}E_{l-k}\left( 0\right) 
\end{equation*}%
where $B_{k,n}^{-1}\left( x\right) $ is the inverse of $B_{k,n}\left(
x\right) $.
\end{corollary}


\begin{thebibliography}{99}
\bibitem{Acikgoz1} M. Acikgoz and S. Araci, On the generating function of
the Bernstein polynomials, \textit{Numerical Analysis and Applied Mathematics%
}, AIP, pp. 1141-1143, \textbf{2010.}

\bibitem{Acikgoz2} M. Acikgoz and S. Araci, New generating function of \
Bernstein type polynomials for two variables, \textit{Numerical Analysis and
Applied Mathematics}, AIP, pp. 1133-1136, \textbf{2010.}

\bibitem{Acikgoz3} M. Acikgoz, and S. Araci, A study on the integral of the
product of several type Bernstein polynomials, \textit{IST Transaction of
Applied Mathematics Modelling and Simulation}, \textbf{2010}, vol. 1, no.
1(2), ISSN 1913-8342, pp. 10-14.

\bibitem{Araci1} S. Araci, Novel identities for $q$-Genocchi numbers and
polynomials, \textit{Journal of Function Spaces and Applications}, Volume
2012 (\textbf{2012}), Article ID 214961, 13 pages.

\bibitem{Araci4} S. Araci, Novel identities involving Genocchi numbers and
polynomials arising from applications from umbral calculus, \textit{Applied
Mathematics and Computation} 233 (\textbf{2014}) 599-607.

\bibitem{Araci5} S. Araci, A. Bagdasaryan, C. \"{O}zel, and H. M.
Srivastava, New Symmetric Identities Involving $q$-Zeta type function,%
\textit{\ Appl. Math. Inf. Sci.} 8, No. 6, 1-6 (\textbf{2014}) (To be
published).

\bibitem{Araci-Acikgoz-Kilicman} S. Araci, M. Acikgoz and A. Kilicman,
Extended $p$-adic $q$-invariant integrals on $%
\mathbb{Z}
_{p}$ associated with applications of umbral calculus, \textit{Advances in
Difference Equations} 2013, \textbf{2013}:96

\bibitem{Araci} S. Araci, M. Acikgoz and E. \c{S}en, On the extended Kim's $p
$-adic $q$-deformed fermionic integrals in the $p$-adic integer ring, 
\textit{Journal of Number Theory} 133 (\textbf{2013}) 3348--3361.

\bibitem{Araci2} S. Araci, M. Acikgoz and H. Jolany, On a $p$-adic
interpolating function associated with modified Dirichlet's type of twisted $%
q$-Euler numbers and polynomials with weight alpha, \textit{Journal of
Classical Analysis} Volume 2, Number 1 (\textbf{2013}), 35--48.

\bibitem{Araci3} E. Cetin, M. Acikgoz, I. N. Cangul and S. Araci, A note on
the ($h,q$)-Zeta-type function with weight $\alpha $, \textit{Journal of
Inequalities and Applications} 2013, 2013:100.

\bibitem{Qi} F. Qi, Integral representations and properties of Stirling
numbers of the first kind, \textit{Journal of Number Theory}, Volume 133,
issue 7 (July, \textbf{2013}), p. 2307-2319

\bibitem{Kim4} T. Kim, $q$-Volkenborn integration, \textit{Russ. J. Math.
Phys.} 9 (2002), no. 3, 288--299.

\bibitem{Kim7} T. Kim, Symmetry $p$-adic invariant integral on $%
\mathbb{Z}
_{p}$ for Bernoulli and Euler polynomials, \textit{J. Difference Equ. Appl.}
14 (2008), no. 12, 1267-1277.

\bibitem{Kim6} T. Kim and D. S. Kim, Applications of umbral calculus
associated with $p$-adic invariants integral on $%
\mathbb{Z}
_{p}$, \textit{Abstract and Applied Analysis,} Vol. 2012 (2012), Article ID
865721, pages 12.

\bibitem{Kim8} T. Kim, D. S. Kim, T. Mansour, S.-H. Rim, M. Schork, Umbral
calculus and Sheffer sequences of polynomials, \textit{J. Math. Phys}. 54 ,
083504 (2013).

\bibitem{Kim5} T. Kim, Some formulae for the $q$-Bernstein polynomials and $%
q $-deformed binomial distributions, \textit{Journal of Computational
Analysis and Applications}, Vol. 14, No.5, 917-933, 2012.

\bibitem{Kim 2} T. Kim, A note on $q$-Bernstein polynomials, \textit{Russian
Journal of Mathematical Physics} \textbf{18} (2011), no. 1.

\bibitem{Kim 5} T. Kim, Analytic continuation of multiple $q$-zeta functions
and their values at negative integers, \textit{Russ. J. Math Phys. }11 (%
\textit{2004}), 71-76.

\bibitem{Kim-Lee-Chae Jang} T. Kim, L-C. Jang and H. Yi, Note on the
modified $q$-Bernstein polynomials, \textit{Discrete Dynamics in Nature and
Society} \textbf{2010 (}2010\textbf{), }Article ID 706483, 12 pages.

\bibitem{Kim 1} T. Kim, J. Choi and Y. H. Kim, On the $k$-dimensional
generalization of $q$-Bernstein polynomials, \textit{Proceedings of the
Jangjeon Mathematical Society}, vol. 14, no. 2, pp. 199-207, 2011.

\bibitem{Kim 3} T. Kim, J. Choi and Y-H. Kim, $q$-Bernstein polynomials
associated with $q$-Stirling numbers and Carlitz's $q$-Bernoulli numbers, 
\textit{Abstract and Applied Analysis}, \textbf{2010} (2010), Article ID
150975, 11 pages.

\bibitem{Kim9} T. Kim, $q$-Euler numbers and polynomials associated with $p$%
-adic $q$-integrals, \textit{Journal of Nonlinear Mathematical Physics, }%
Volume 14, Number 1 (\textbf{2007}), 15-27.

\bibitem{Would} H. W. Gould, Explicit formulas for Bernoulli numbers, 
\textit{Amer. Math. Monthly} 79(1972), 44-51.

\bibitem{Choi 1} J. Choi and T. Y. Seo, Identities involving series of the
Riemann Zeta function, \textit{Indian J. Pure Appl. Math.} 30 (1999),
649-652.

\bibitem{Choi 2} J. Choi and H. M. Srivastava, Certain classes of series
associated with the Zeta function and multiple Gamma functions, \textit{J.
Comput. Appl. Math.} 118 (2000), 87-109.

\bibitem{Choi 3} J. Choi, Remark on the Hurwitz-Lerch zeta function, \textit{%
Fixed Point Theory and Applications} 2013, 2013:70.

\bibitem{Choi 4} J. Choi and H. M. Srivastava, Integral Representations for
the Gamma function, the Beta Function, and the Double Gamma Function, 
\textit{Integral Transforms Spec. Funct.} 20(11) (2009), 859--869.

\bibitem{Srivastava} H. M. Srivastava, Some generalizations and basic (or $q$%
-) extensions of the Bernoulli, Euler and Genocchi polynomials, \textit{%
Appl. Math. Inf. Sci.} 5, 390-444 (2011).

\bibitem{Srivastava1} H. M. Srivastava, Some formulas for the Bernoulli and
Euler polynomials at rational arguments, \textit{Math. Proc. Cambridge
Philos. Soc.} 129 (2000), 77--84.

\bibitem{Srivastava2} H. M. Srivastava and J. Choi, \textit{Zeta and }$q$%
\textit{-Zeta Functions and Associated Series and Integrals, }Elsevier
Science Publishers, Amsterdam, London and New York, 2012.

\bibitem{Kim 4} V. Gupta, T. Kim J. Choi, and Y. H. Kim, Generating
functions for $q$-Bernstein, $q$-Meyer-K\"{o}nig-Zeller and $q$-Beta basis, 
\textit{Automation Computers Applied Mathematics} Vol. 19 (2010) No. 1, pp.
119-122.

\bibitem{Kellner} B. C. Kellner, On irregular prime power divisors of the
Bernoulli numbers, \textit{Mathematics of Computation}, Volume 76, Number
257, January 2007, Pages 405--441.

\bibitem{Kruchinin} D.V. Kruchinin and V. V. Kruchinin, Application of a
composition of generating functions for obtaining explicit formulas of
polynomials, \textit{Journal of Mathematical Analysis and Applications},
Volume 404, Issue 1, 1 August \textbf{2013}, Pages 161--171.

\bibitem{Mezo} I. Mez\H{o} and A. Dil, Hyperharmonic series involving
Hurwitz zeta function, \textit{J. Number Theory.} 130(2) (2010), 360-369.

\bibitem{Buyukyazici1} I. Buyukyazici and E. Ibikli, Bernstein polynomials
of two variable functions, \textit{Graduate School of Natural and Applied
Sciences}, Department of Mathematics, 1999, 49 pages, Ankara, Turkey.

\bibitem{Buyukyazici2} I. Buyukyazici and E. Ibikli, The approximation
properties of generalized Bernstein polynomials of two variables, \textit{%
Applied Math. and Comput.} \textbf{156} (2004) 367-380.

\bibitem{Oruc2} H. Oruc and G. M. Phillips, A generalization of the
Bernstein polynomials, \textit{Proceedings of the Edinburgh Mathematical
society} (1999) 42, 403-413.

\bibitem{Ostrovska} S. Ostrovska, On the $q$-Bernstein polynomials, \textit{%
Adv. Stud. Contemp. Math.} 11 (2) (2005), 193-204.

\bibitem{Phillips} G. M. Phillips, A survey of results on the $q$-Bernstein
polynomials, \textit{IMA Journal of Numerical Analysis} Advance Access
published online on June 23, (2009), 1-12, doi:10.1093/imanum/drn088.

\bibitem{SimsekAcikgoz1} Y. Simsek and M. Acikgoz, A new generating function
of $q$-Bernstein-type polynomials and their interpolation function, \textit{%
Abstract and Applied Analysis}, volume 2010, Article ID 769095, 12 pages,
doi: 10.1155/2010/769095.01-313.
\end{thebibliography}
\end{document}